\documentclass[11pt]{amsart}
\usepackage{amsmath,amsthm,amsfonts,amssymb,amscd}
\usepackage[latin1]{inputenc}
\usepackage{palatino}

\usepackage{graphicx}

\newtheorem{Theorem}{Theorem}[section]

\newtheorem{Proposition}[Theorem]{Proposition}
\newtheorem{Lemma}[Theorem]{Lemma}
\newtheorem{Claim}[Theorem]{Claim}
\theoremstyle{Definition}
\newtheorem{Definition}[Theorem]{Definition}
\newtheorem{Example}[Theorem]{Example}

\theoremstyle{Remark}
\newtheorem{Remark}[Theorem]{Remark}

\def\leaderfill{\leaders\hbox to .8em{\hss .\hss}\hfill}
\def\_#1{{\lower 0.7ex\hbox{}}_{#1}}

\def\fa{{\mathcal{F}}}

\def\vr{{\varphi}}

\def\sing{\operatorname{{sing}}}

\def\cod{\operatorname{{cod}}}
\def\codim{\operatorname{{cod}}}

\title{Integrable deformations of local analytic fibrations with singularities}
\author{Dominique Cerveau and Bruno Scárdua}

\begin{document}

\maketitle

\begin{abstract}
We study analytic integrable deformations of the germ of a holomorphic foliation given by $df=0$ at the origin $0 \in \mathbb C^n, n \geq 3$. We consider the case where $f$ is a germ of an irreducible and reduced holomorphic function. Our central hypotheses is that, {\em outside of a dimension $\leq n-3$ analytic subset $Y\subset X$,   the analytic hypersurface $X_f : (f=0)$ has only normal crossings singularities}. We then prove that, as germs,  such deformations also exhibit a holomorphic first integral, depending analytically on the parameter of the deformation. This applies to the study of integrable germs writing as $\omega = df + f \eta$ where $f$ is quasi-homogeneous. Under the same hypotheses for $X_f : (f=0)$ we prove that $\omega$ also admits a holomorphic first integral. Finally, we conclude that an integrable germ $\omega = adf + f \eta$ admits a holomorphic first integral provided that:
(i) $X_f: (f=0)$ is irreducible with an isolated singularity at the origin $0 \in \mathbb C^n, n \geq 3$; \,
(ii) the algebraic multiplicities of $\omega$ and $f$ at the origin satisfy $\nu(\omega) = \nu (df)$. In the case of an isolated singularity for $(f=0)$ the writing $\omega = adf + f \eta$ is always assured so that we conclude  the existence of a holomorphic first integral.  Some questions related to Relative Cohomology are naturally considered and not all of them answered.

\end{abstract}


\section{Introduction and main results}
\label{section:introduction}
The problem of integrability of differential equations in the real  context goes back to H.Poincaré and Dulac. In the analytic case it is natural to consider the complexification of the equation and then we are just one step away from the holomorphic foliations framework. These are objects that can be described by integrable systems of (holomorphic) one-forms. Under this viewpoint probably the most important result is Malgrange's work (\cite{malgrangeI,malgrangeII}, relating the dimension of the singular set of the system with the existence of a holomorphic first integral for it. This is one of the motivations for this work.

\subsection{Main results}

We consider $f\colon \mathbb C^n , 0 \to \mathbb C,0$ a germ of a holomorphic function at the origin $0\in \mathbb C^n, n \geq 3$. The corresponding germ of an  analytic hypersurface $(f=0)$ is denoted by $X_f$.  The singular set of the hypersurface $X_f$ will be denoted by $\sing(X_f)$. We will denote by $X_f^*= X_f - \sing (X_f)$ the smooth part of $X$. Next we give a pleonastic definition of our main hypothesis:

\begin{Definition}
{\rm We shall say that $X_f$ {\it has only ordinary singularities off a codimension $\geq 3$ subset} if  there exists an analytic subset $(Y,0)\subset (X_f,0)$ of dimension at most $n-3$, such that outside of $Y$ the only singularities of $(X_f,0)$ are normal crossings.
}
\end{Definition}

We will assume that $f$ is {\it reduced} (if $g \in \mathcal O_n$ is such that $g \big|_{X_f} \equiv 0$ then $f\big| g$ in $\mathcal O_n$.). In this case the singular set of $X_f$ is given by $\sing(X_f) = \sing(f) = \{p \in (\mathbb C^n,0): df(p)=0\}$. Indeed, it is well-known (\cite{Milnor}) that the singular points of $f$, i.e.,  the zeroes of $df$, are contained in the fiber $f^{-1}(0)$. We consider the germ of an integrable one-form $\omega \in \Omega^1(\mathbb C^n, 0)$. Then $\omega =0$ defines a codimension-one holomorphic foliation $\fa(\omega)$ germ at $0 \in \mathbb C^n$. The hypersurface $X_f$ is $\fa(\omega)$-invariant if, and only if, $\omega \wedge df/f$ is holomorphic (\cite{cerveau-berthier,saito}). This is the case of integrable one-forms that write as
\[
\omega = a df + f \eta
\]
with $a \in \mathcal O_n$ and $\eta \in \Omega^1(\mathbb C^n,0)$.

In a certain sense this is the most natural writing for a one-form $\omega$ that leaves $X_f$ invariant (see \S~\ref{section:isolatedsingularity} and \S~\ref{section:relativecohomology}).

For $\eta$ small enough (in the sense of Krull topology \cite{kaup-kaup}, \cite{gunning-rossi}) and $a \in \mathcal O_n ^*$ unit, we may see $\fa(\omega$ as an integrable deformation of the holomorphic "fibration" $\fa(df)$, given by $f=const.$. If for instance $f$ has an isolated singularity at $0 \in \mathbb C^n, n \geq 3$, then any $\omega$ that leaves $X_f : (f=0)$ invariant must write as above, $\omega = adf + f \eta$.
In particular, $\omega$ may come from an analytic deformation of $\omega_0=df$, under some geometrical condition (e.g., if $\nu(\omega= \nu(df)$ as explained in the text itself).

We will consider the following situation:
$\{\omega_t, t\in (\mathbb C,0)\}$ is an analytic deformation of $\omega_0 = df$ such that each one-form $\omega_t\in \Omega^1 (\mathbb C^n,0)$ is integrable, $\omega_t \wedge d \omega_t=0$.

We then prove in \S~\ref{section:Lesaito}:
\begin{Theorem}
\label{Theorem:1}
Assume that the germ $f\in \mathcal O_n, n \geq 3$ is reduced, $X_f$ is irreducible and  has only normal crossings  singularities off a codimension $\geq 3$ subset. Let $\{\omega_t\}_{t \in \mathbb C,0}$ be an analytic deformation of $\omega_0=df$ at $0 \in \mathbb C^n$. Then for any $t \in (\mathbb C,0)$ close enough to $0$, the one-form $\omega_t$ admits a holomorphic first integral. Indeed, there is a germ of a holomorphic function $F \colon \mathbb (\mathbb C^n \times \mathbb C , 0) \to (\mathbb C,0), \, (x,t) \mapsto F(x,t)$ such that:
\begin{enumerate}
\item $F(x,0)=f(x)$
\item $F_t \colon x \mapsto F(x,t)$ is a first integral for $\omega_t$.
\end{enumerate}
\end{Theorem}

\begin{Remark}
{\rm Regarding the hypotheses in Theorem~\ref{Theorem:1}, we observe that:
\begin{enumerate}
 \item {\em It is enough to assume that they hold for the restriction of $\omega$ to a three-dimensional plane, in {\it general position with respect to} $\omega$, in the same sense of {\rm\cite{mattei-moussu}}.}

     Indeed, according to
\cite{mattei-moussu} if such a restriction admits a holomorphic first integral, then the same holds for the form $\omega$. Actually, it is proved that the first integral for the $3$-dimensional plane section admits an extension to a first integral for the form $\omega$ in a neighborhood of the origin $0 \in \mathbb C^n$.
\item In particular, {\em if there is a three-dimension plane section  of $X$ which has an isolated singularity at the origin, then the conclusion of Theorem~\ref{Theorem:1} is valid.}

\end{enumerate}
}
\end{Remark}
As  a corollary we have (see \S~\ref{section:isolatedsingularity}):
\begin{Theorem}
\label{Theorem:2}
Let $f\in \mathcal O_n, n \geq 3$ be a strictly quasi-homogeneous reduced function. Assume that $X_f : (f=0)$  is irreducible and has  only normal crossings singularities off a codimension $\geq 3$ subset. Then, any holomorphic integrable one-form germ $\omega\in \Omega^1(\mathbb C^n,0)$  of the form $\omega = df + f \eta$, with $\eta\in \Omega^1(\mathbb C^n,0)$,  admits a holomorphic first integral.
\end{Theorem}

\begin{Remark}
{\rm A function $f(x_1,...,x_n)$ is {\it quasi-homogeneous} if there exist $d, d_1,...,d_n$ all non-negative such that $f(t^{d_1} x_1,...,t^{d_n} x_n)=t^df(x_1,...,x_n)$. It is {\it strictly} quasi-homogeneous if $d>0$ and $d_j>0, \forall j$.
}

\end{Remark}

Then, as a consequence of our approach and some relative cohomology results based on \cite{saito} and \cite{cerveau-berthier} we obtain, for the case of an isolated singularity (cf. \S~\ref{section:isolatedsingularity}):

\begin{Theorem}
\label{Theorem:3}
 Let  $f \in \mathcal O_n$ be reduced and irreducible, with an isolated singularity at the origin. Let $\omega \in \Omega^1(\mathbb C^n,0), n \geq 3$ be an integrable one-form having  $(f=0)$ as the only invariant hypersurface. Then $\omega$ admits  a (germ of a) holomorphic first integral if, and only if,  $\nu(\omega) = \nu(df)$ at $0$.

\end{Theorem}

In the above statement, $\nu(.)$ stands for the algebraic multiplicity at the origin.

\section{Germs of hypersurfaces with normal crossings}
\label{section:Lesaito}
We consider $(X_f,0)\subset (\mathbb C^{n},0)$ a germ of reduced analytic surface in $(\mathbb C^{n},0)$, defined by $f=0$ where $f \in \mathcal O_{\mathbb C^{n},0}$ is a germ of analytic function in $0 \in \mathbb C^{n}$.
If $n=3$ and   $(X_f,0)$ only has normal crossings singularities off the origin $0 \in \mathbb C^{3}$, then the local fundamental group  of the complement of $(X_f,0)$ in $(\mathbb C^{n},0)$ is abelian. This is a particular case of the more general statement below:

\begin{Theorem}[ L\^e-Saito, \cite{Le-Saito} Main Theorem page 1]
\label{Theorem:LeSaito}
Let $n \geq 3$. Assume that outside of an analytic subset $(Y,0)\subset (X_f,0)$ of dimension at most $n-3$, the only singularities of $(X_f,0)$ are normal crossings. Then the  local fundamental group of the complement of $(X_f,0)$ in $(\mathbb C^{n},0)$ is abelian. The Milnor fiber of $f$ has a fundamental group which is free abelian of rank the number of analytic components of $X_f$ at $0$, minus one. Finally, if $X_f$ is irreducible, then the  fiber $f^{-1}(c), c \ne 0$  is simply-connected.
\end{Theorem}

Throughout this section,  we will consider the case where $X=X_{f_0}$ is irreducible,  given by $f_0=0$ as above  with normal crossings singularities off a codimension $\geq 3$ subset. We study analytic integrable deformations of the one-form $\omega_0=df_0$. Such a deformation writes as
\[
\omega_t = \omega_0 + t \omega_1 + \ldots + t^k \omega_k + \ldots
\]

where $t \in \mathbb C,0$ and the $\omega_j$ are holomorphic in some small neighborhood $U$ of $ 0 \in \mathbb C^n, n \geq 3, \forall j \geq 0$.
The integrability condition $\omega_t \wedge d \omega_t=0$ gives:

\[
\omega_0 \wedge d \omega_0=0
\]
\[
\omega_0 \wedge d\omega_1 + \omega_1 \wedge d\omega_0 =0
\]
\[
\omega_2 \wedge d \omega_0 + \omega_1 \wedge d \omega_1 + \omega_0 \wedge d\omega_2=0
\]
\[
\vdots
\]
In our case $\omega_0 = df_0$, i.e., we have
$d\omega_0=0$ and then
\[
df_0 \wedge d \omega_1=0
\]
and
\[
\omega_1 \wedge d\omega_1 + df_0 \wedge d \omega_2 =0
\]

Now, $\omega_1$ is not necessarily integrable, but we have the following Relative Cohomology Lemma:

\begin{Lemma}
\label{Lemma:1}
Under the above hypotheses we have:
\[
d\omega _1 \wedge df_0 =0 \implies \omega_1 =df_1 + a_1 df_0
\]
for some holomorphic functions $f_1,a_1$ in $U$.

\end{Lemma}

The proof of Lemma~\ref{Lemma:1} is given in \S~\ref{section:relativecohomology} as a consequence of Theorem~\ref{Theorem:LeSaito} and some relative cohomology techniques  adapted from \cite{cerveau-berthier}.



\begin{Remark}

The writing $\omega_1 = df_1 + a_1 df_0$ is not unique but if $\omega_1=
d\tilde f_1 + \tilde a_1 df_0$ then $\tilde f_1 = f_1 + \vr(f_0)$ and $a_1 = \tilde a_1 + \vr^\prime (f_0)$ for some one-variable holomorphic germ $\vr(z)$.

\end{Remark}

Returning to the deformations we obtain
\[
\omega_t= \omega_0 + t \omega_1 + t^2 \omega_2 + \ldots = df_0 + t(df_1 + a_1 df_0) + t^2 \omega_2 + \ldots =
(1+ t a_1) df_0 + t df_1 + t^2 \omega_2 + \ldots
\]
Therefore
\[
\frac{1}{1+ t a_1} \omega_t = df_0 + \frac{t}{1+ t a_1} df_1 + \frac{t^2}{1+ t a_1} \omega_2 + \ldots
=d(f_0 + t f_1)  + t^2 \tilde \omega_ 2 + \ldots
\]

We put
\[
\tilde \omega_t = \frac{1}{1+ t a_1} \omega_t
\]
recalling that $ 1 + t a_1$ is a unit and $\tilde \omega_t$ is also integrable holomorphic. Then we write
\[
\tilde \omega_t = df_0 + tdf_1 + t^2 \tilde \omega_2 + \tilde \omega_3 + \ldots,
\]
\[
d\tilde \omega_ t = t^2 d \tilde \omega_2 + t ^3 d \tilde \omega_3 + \ldots
\]
From the integrability condition we then obtain
$df_0 \wedge d \tilde \omega_2=0$. As above we obtain
\[
\tilde \omega_ 2 = df_2 + a_2 df_0
\]
for some holomorphic $a_2, f_2 \colon U \to \mathbb C$.
Then
\[
\tilde \omega_ t = df_0 + t df_1 + t^2 (df_2 + a_2 df_0) + t^3 \tilde \omega_3 +\ldots =
(1+ t^2 a_2) df_0 + tdf_1 + t^2 df_2 + t^3 \tilde \omega_ 3 + \ldots
\]
Hence, as above, we can define holomorphic integrable
\[
{\tilde {\tilde{\omega_t}}} := \frac{1}{1+ t^2 a_2} \tilde \omega_t =
df_0 + tdf_1  + t^2 df_2 t^3 \tilde {\tilde {\omega_3}} + \ldots
\]

Inductively proceeding like this we obtain a formal unit $\hat G$ such that
\[
\frac{1}{\hat G}\omega_t = df_0 + \sum\limits_{j=1} ^\infty t^j df _j
\]
for some holomorphic functions $f_j \colon U \to \mathbb C, j \geq 1$.

Thus we can write $\omega_ t=  \hat G(x,t). d_x \hat F(x,t)$ in the obvious sense for formal function $\hat F(x,t)$ (indeed, $\hat F$ is what is referred to as "transversely formal" in the sense of Mattei-Moussu).

Remark that there is a formal series $\hat h(x,t)$ such that $\omega_t + \hat h dt= \hat G d_{(x,t)} \hat F$.

Now we consider the pair $\{\Omega, dt\}$ where $\Omega(x,t)=\omega_t(x)$, defined in $(U\times \mathbb C,0)\subset (\mathbb C^{n+1},0)$. We claim that this is an integrable system: Indeed, $d_{(x,t)}\Omega= d_x \omega_t + \frac{\partial \omega_t}{\partial t} dt $ so that
\[
\Omega \wedge d \Omega \wedge dt = \omega_t \wedge d_x \omega_t \wedge dt +
\omega_t \wedge \frac{\partial \omega_t}{\partial t} dt \wedge dt= 0
\]
because we have $\omega_t \wedge d_x\omega_t=0$.

We also claim that
$\{\omega, dt\} = \{ d\hat F, dt\}$ at the level of formal modules. This is immediate from the expression $\Omega = \hat G d\hat F - \hat h dt$ and from the fact that $\hat G$ is a unit. Now we recall the following result of Malgrange:

\begin{Theorem}[Malgrange, \cite{malgrangeII}]
Let be given $p$ germs of holomorphic one-forms \linebreak $\omega_1,...,\omega_p\in \Omega^1 (\mathbb C^n,0)$ at the origin $0 \in \mathbb C^n$ with $ 1 \leq p \leq n$. Denote by $\mathcal S$ the germ of analytic set given by the zeros of $\omega_1 \wedge \ldots \omega_p$. Assume that one of the conditions below is verified:

\begin{itemize}
\item[{\rm(i)}] $\cod \mathcal S \geq 3$ and $d\omega_j \wedge \omega_1 \ldots \wedge\omega_p=0, \forall j$.
    \item[{\rm (ii)}] $\cod S \geq 2$ and $\{\omega_1,...,\omega_p\}$ is {\it formally integrable}  (i.e., there exist $\hat f_i, \hat g_{ij} \in \hat {\mathcal O}_n$ such that $\omega_i = \sum\limits_j \hat g_{ij} d\hat f_j, \forall i$ and $\det(\hat g_{ij} (0) \ne 0$).

        \end{itemize}
        Then $\{\omega_1,...,\omega_p\}$ is {\it integrable} meaning that there exist $ f_i, g_{ij} \in \mathcal O_n$ (convergent) such that $\omega_i = \sum\limits_j g_{ij} df_j, \forall i$ and $\det(g_{ij}(0)) \ne 0$.
        \end{Theorem}

        We shall apply this result to the system $\{\Omega, dt\}$. Recall that $\hat {\mathcal O}\{\Omega, dt\}= \{ d \hat F, dt\}$, that is, $\{\Omega, dt\}$ is formally integrable. Notice that $\Omega \wedge dt = d_x f_o \wedge dt + \sum\limits_{j \geq 1} t^j \omega_j \wedge dt$. Since $\codim \sing (df_0) \geq 2$ we conclude that $\cod \sing (\Omega \wedge dt) \geq 2$.

Thus, by (ii) in Malgrange's theorem above we have that $\{d \hat F, dt\}$ admits a holomorphic (convergent) pair $\{dF, dt\}$. i.e., there is a holomorphic function $F(x,t)$ in a neighborhood of $0 \in \mathbb C^n$ such that  $\Omega= a. d_x F(\cdot, t) + b dt$ for some holomorphic functions $a,b$ where $a$ is a unit.  Therefore $\omega_t = a .d_x F(x,t)$. This completes the proof of Theorem~\ref{Theorem:1}.

\qed
\section{Necessity of hypothesis, proof of Theorem~\ref{Theorem:2}}
\label{section:example}

We shall now discuss the necessity of our central hypothesis, about the
normal crossings for the singularities of the hypersurface $X : (f=0)$.
\begin{Example}
\label{Example:sharp}
{\rm
We consider $f_0 = x^3 + y^2$ and the corresponding cusp $f_0=0$ in $\mathbb C^2,0$.

The cylinder generated by this cusp in $\mathbb C^3,0$ is a hypersurface $X$ with singular set irreducible of codimension two. We consider deformations of the form
\[
\omega_t = d(y^2 + x^3) + tx (2xdy - 3y dx)
\]

From \cite{Cerveau-Mattei}, \cite{Rafik-Loray} it is known that there are no holomorphic first integrals for such a generic deformation. Embedding this on $\mathbb C^3,0$ we conclude that the hypothesis of normal crossings for the hypersurface singularities cannot be dropped.

}
\end{Example}

 \begin{proof}[Proof of Theorem~\ref{Theorem:2}] We prove Theorem~\ref{Theorem:2} as an application of Theorem~\ref{Theorem:1} for the case of quasi-homogeneous hypersurface $X_f$.
 We start with a holomorphic integrable germ of one-form at $0 \in \mathbb C^n, n \geq 3$ of the following type
\[
\omega = df + f \eta
\]
where $f$ is a strictly quasi-homogeneous function with normal crossings. This means that:
\begin{enumerate}
\item $\exists\,  d,d_1,...,d_n > 0,$ such that $f(t^{d_1} x_1,...,t^{d_n} x_n) = t^d f(x_1,...,x_n);$

    \item $X_f$ is an analytic irreducible hypersurface with ordinary (normal) crossings singularities off a codimension $\geq 3$ subset.
        \end{enumerate}

        Let us see how to consider (embed) $\omega$ as (into) a deformation of $\omega_0 = df$.

        Take the map $\sigma_t \colon \mathbb C^n,0 \to \mathbb C^n,0$ defined for $ t \ne 0$ by

        $\sigma _ t (x) = (t^{d_1} x_1,...,t^{d_n} x_n)$. Then we have
        $\sigma _t ^* f = f \circ \sigma_t = t^d f$ and
        \[
        \sigma_t ^* \omega = \sigma_t ^* ( df+ f \eta) =
        t^d df + t^d f \sigma_t ^* (\eta).
        \]

        We then define $\omega_t := \frac{1}{t^d} \sigma_t ^* \omega=
        df + f \sigma_t ^* (\eta)$.

        Notice that, because $d_j >0, \forall j$ we have $\sigma_t ^* (\eta) = t \tilde \eta$ for some holomorphic one-form $\tilde \eta (x,t)$.Therefore $\omega _t = \frac{1}{t^d} \sigma_t ^* \omega =
        df + t f \tilde \eta$. Thanks to this form, $\omega_t$ is an analytic deformation of $\omega_0= df$, such that $\omega_1 = \omega$. Theorem~\ref{Theorem:2} then follows from the fact that $\sigma_t$ defines an analytic diffeomorphism taking $\fa(\omega_1)$ onto $\fa(\omega_t)$ for all $t \ne 0$, and from Theorem~\ref{Theorem:1} the foliation $\fa(\omega_t)$ admits a holomorphic first integral for $t$ close enough to $0$.
        \end{proof}

\section{The case of an isolated singularity}
\label{section:isolatedsingularity}

We consider an integrable germ of holomorphic one-form $\omega \in \Omega^1(\mathbb C^n,0)$ with an invariant hypersurface $X_f: (f=0)$ such that
\begin{itemize}
\item[{\rm (Is.1)}] $\cod \sing (\omega) \geq 2$
\item[{\rm (Is.2)}] $X_f$ has an isolated singularity at $0\in \mathbb C^n, n \geq 3$.
\end{itemize}

Since $(f=0)$ is $\omega$-invariant we have that $f \big| \omega \wedge df$in $\Omega^2 (\mathbb C^n,0)$, i.e., $\frac{1}{f} \omega \wedge df$ is holomorphic. Because of the above hypotheses (i) and (ii) we can indeed write $ \omega = adf + f \eta$ for some holomorphic function germ $a$ and holomorphic one-form germ $\eta$. This is the content of the following lemma:

\begin{Lemma}
If $X: (f=0)$ is irreducible, reduced and has an isolated singularity (at the origin $0 \in \mathbb C^n, n \geq 3$) then $\omega$ can be written $(*) \, \, \omega = adf + f \eta$
for some holomorphic $a, \eta$.

\end{Lemma}
\begin{proof}
This is related to the Dolbeault Cohomology of $\mathbb C^n\setminus 0$ and the corresponding vanishing theorem of Cartan: $H^1(\mathbb C^n- \{0\},\mathcal O)=0$ if $ n \geq 3$. Locally at any point off the origin we may write $\omega$ as in (*). This gives an open cover $\bigcup\limits_{j \in \mathbb N} U_j$ of a punctured neighborhood
$U^* = U \setminus 0$ of $ 0 \in \mathbb C^n$. We may assume that $U$ is a polydisc centered at the origin. For each open set $U_j \subset U$, the restriction $\omega\big|_{U_j}$ writes $\omega= a_j df + f \eta_j$ for some holomorphic $a_j,\eta_j$ in $U_j$. We can assume that each intersection $U_i \cap U_j \ne \emptyset$ is connected. If $U_i \cap U_j\ne\emptyset$ then on $U_i \cap U_j$ we have
$a_ i df + f \eta_i = a_j + f \eta_j$. Thus we have $(a_i - a_j ) df = f (\eta_j - \eta_i)$. So, because $df$ does not vanish on $U_i \cap U_j\subset U \setminus 0$ we have $f\big| a_i - a_j$ in $\mathcal O (U_i \cap U_j)$. We write then $a_i  - a_j = f h_{ij}$ for some $h_{ij} \in \mathcal O(U_i \cap U_j)$. The data $\{h_{ij}, U_i\}$ defines an additive cocycle in $U^*$, so that by Cartan's theorem this cocycle has a solution (\cite{gunning1,gunning-rossi}). This means that there exist $h_j \in \mathcal O (U_j)$ such that on each non-empty intersection $U_i \cap U_j$ we have $h_{ij} = h_i - h_j$. Therefore $a_i - fh_i = a_j - f h_j$. Thus we may define $a\in \mathcal O(U^*)$ by $a\big|_{U_j} = a_j - f h_j$. From the above equations we have
that $h_i df + \eta _i = h_j df + \eta _j$. We then define $\eta$ in $U^*$ by setting $\eta\big|_{U_j} = h_j df + \eta_j$. Then on each $U_j$ we have
$\omega = a_j df + f \eta_j = adf + \eta$. By classical Hartogs' extension theorem (\cite{gunning1,gunning2}) the function $a$ and the form $\eta$ extend to $U$ and we write $\omega = adf + f \eta$ in $U$.

\end{proof}

\begin{Remark}
{\rm The above lemma may also follow from the following argumentation, based on Saito-De Rham division lemma (\cite{saito}):
From the invariance of $(f=0)$, where $f$ is reduced, we have
$\omega \wedge df = f \theta$ for some holomorphic $\theta \in \Omega^2(\mathbb C ^n,0)$. This means that $\omega \wedge df=0$ in the quotient ring   $\mathcal O_n /f$. Then, because $\sing (f)=\{0\}\subset \mathbb C^n$ and $ n \geq 3$, we have from \cite{saito} (page 166)  that $\omega = a df $ in $\mathcal O_n /f$. In other words, $\omega = adf + f\eta$ for some holomorphic $\eta \in \Omega^1(\mathbb C^n,0)$.
}
\end{Remark}

Now we proceed under the hypothesis that $\omega = adf + \eta$. For the case of an isolated singularity this is always the case. Indeed, we can say more:
\begin{Claim}
Assume that $f$ has an isolated singularity at the origin and that:
 \begin{itemize}
 \item[{\rm (Is.3)}] the algebraic multiplicities of $\omega$ and $f$ at the origin satisfy $\nu(\omega) = \nu(df)$.
      \end{itemize}

      Then the function $a$ is a unit.
\end{Claim}
\begin{proof}
We have $\omega = adf + f \eta$. From this equation we conclude that $a$ is a unit, simply by comparing the orders of $f$ and $df$ at the origin, plus using the fact that $\omega$ and $df$ have the same order at the origin.
\end{proof}

From now on we assume that $\omega = a df + f\eta$ where the function $a$ is a unit. Next we show that  $\sing(\omega)=\{0\}$. We can suppose that $a=1$; note that $f$ is a submersion outside $\{0\}$. This implies that $\sing (\omega)\cap f^{-1}(0)\subset\{0\}$. Suppose that $\sing(\omega)$  contains a curve parametrized by  $t \mapsto \gamma(t), t \in (\mathbb C,0)$.  Then $f\circ \gamma(t) \not \equiv 0$ and up to reparametrization we can suppose that  $f(\gamma(t))=t^p$ for some $0< p \in \mathbb N$.
Then $0=\gamma^*(\omega) =  p t^{p-1} dt + t^p \gamma^*(\eta)$. This implies $0=p dt + t \gamma^*(\eta)$,  a contradiction.

\begin{Remark}
\label{Remark:equivalentIs.}
{\rm Now we examine  another condition:
\begin{itemize}
\item[{\rm (Is.3')}]For a generic plane section $E\colon (\mathbb C^2,0) \hookrightarrow (\mathbb C^n, 0)$, the restriction $E^*(\omega)\in \Omega^1(\mathbb C^2,0)$ defines a foliation which is non-dicritical with a singularity of generalized curve (\cite{C-LN-S1}) type at the origin $0 \in \mathbb E^2$.
    \end{itemize}

Then, from  \cite{C-LN-S1} we have that (Is.3') $\implies$ (Is.3).
Indeed, it is possible to give some further conditions on $\fa(\omega)$ in order to conclude that conditions (Is.3) and (Is.3') are equivalent (see the paragraph preceding Theorem~\ref{Theorem:3} below).
}
\end{Remark}

\subsection{Conclusions}
Let us collect our conclusions from the previous discussion:
\subsubsection{Isolated singularity}
For the case of an irreducible  hypersurface $X: (f=0)$ at $\mathbb C^n,0$ ( $ n \geq 3)$ with an  isolated singularity at the origin, we obtain:

\begin{Proposition}
\label{Proposition:preliminary} Given a holomorphic integrable one-form  $\omega= adf + f \eta$ at $0 \in \mathbb C^n, n \geq 3$ with $\cod \sing (\omega) \geq 2$. Then we have $\sing(\omega) = \sing(df) \subset \{0\}$ and we have a germ of a non-constant holomorphic first integral for $\omega$. \end{Proposition}

In particular, let be given a holomorphic integrable germ of a one-form $\omega= adf + f \eta$ with $a\in \mathcal O_n$. Assume that $X_f$ is the only invariant hypersurface and that for a generic plane section
$E\colon (\mathbb C^2,0)\hookrightarrow (\mathbb C^n,0)$, the induced foliation $E^*\fa(\omega)$ is a non-dicritical generalized curve in the sense of \cite{C-LN-S1}. Then the algebraic multiplicities of $\omega$ and $f$ at the origin satisfy $\nu(\omega) = \nu(df)$. Therefore $a$ is a unit and we drop on the preceding case, i.e., we may assume that $\omega= df + f \tilde \eta$. In particular, we conclude:

\begin{Proposition}
\label{Proposition:isolatedsing}
Given a holomorphic integrable one-form  $\omega$ at $0 \in \mathbb C^n, n \geq 3$, assume that:
 \begin{enumerate}
 \item $X_f: (f=0)$ is invariant, where $f$  has an isolated singularity at the origin.
 \item The algebraic multiplicities of $\omega$ and $f$ at the origin satisfy $\nu(\omega) = \nu(df)$.

 \end{enumerate}
 Then there exists a  germ of a holomorphic first integral for $\omega$. \end{Proposition}
\begin{Remark}
{\rm In dimension $n \geq 3$ the fact that $f$ has an isolated singularity at the origin already implies that $X_f$ is irreducible. Condition (2) above is satisfied if $X_f$ is the only invariant hypersurface and if  generic plane sections of $\fa(\omega)$ are non-dicritical generalized curves.
}
\end{Remark}

\begin{proof}[Proof of Theorem~\ref{Theorem:3}]

 In view of Proposition~\ref{Proposition:isolatedsing} (see also Remark~\ref{Remark:equivalentIs.}) it remains to show the "only if" part.
Assume that  $\omega$ has a holomorphic first integral and $(f=0)$ is the only invariant hypersurface. Let us prove that  $\nu(\omega)= \nu(df)$. Indeed,  if $F$ is a holomorphic first integral for $\omega$ then we can assume that $F=g.f$ for some holomorphic $g$.
Since $F$ is a first integral and $\cod \sing (dF) \geq 2$, we can write (\cite{saito}) \, $\omega = b. dF$ for some unit $b$. Taking derivatives  we
obtain $\omega = b .gdf + fb.dg$. Notice that $b.g$ is also a unit. Therefore, we must have $\nu(\omega) = \nu(df)$.

\end{proof}

\section{Some relative cohomology}
\label{section:relativecohomology}
Let us give now in  details  the proof of Lemma~\ref{Lemma:1}.
We consider $f \in \mathcal O_n, n \geq 3$ with $f(0)=0$ and put $X_f: (f=0)\subset (\mathbb C^n,0)$. We consider $\omega \in \Omega^1 (\mathbb C^n,0)$ an integrable germ of holomorphic one-form. We assume that $X_f$ is irreducible.

\begin{Lemma}
\label{Lemma:poincare}
Assume that $\cod \sing(df) \geq 2$. The following conditions are equivalent:
\begin{enumerate}
\item $d \omega \wedge df =0$.
\item $\omega$ is closed on each fiber of $f$.
\end{enumerate}

\end{Lemma}
\begin{proof}
It is sufficient to prove the lemma at a generic point (for ad-hoc representatives of our germs). By Poincar\'e lemma, Lemma~\ref{Lemma:poincare} is true for a submersion. In particular, it is true at a generic point for $f$.
\end{proof}

We recall that if $f$ is reduced then $\cod \sing(df) \geq 2$.
\begin{Proposition}
\label{Proposition:integration}
Assume that $\cod \sing(f) \geq 2$.
Then the following conditions are equivalent for $\omega \in \Omega^1 (\mathbb C^n,0)$:
\begin{itemize}
\item[{\rm(i)}] $ d\omega \wedge df=0$ and $\oint _ \gamma \omega =0$ for each cycle
$\gamma$ contained in a non-singular fiber $f^{-1}(c), c \ne 0$.

\item[{\rm(ii)}] $\omega = adf + dh$ for some $a, h \in \mathcal O_n$.

\end{itemize}

\end{Proposition}
\begin{proof}

Since the sense (ii) $\implies$ (i) is trivial, we shall assume that $d \omega \wedge df =0$ and $\oint _\gamma \omega =0$ for all cycle $\gamma$ contained in the non-singular fibers of $f$. By the preceding lemma we have that the restriction of $\omega$ to these fibers is closed.

 \begin{Claim}
 There exist  holomorphic functions $a,h \colon \mathbb C^n \setminus X_f , 0 \to \mathbb C$ such that $\omega = adf + dh$.
 \end{Claim}

\begin{proof}
We will follow an integration argument like in the proof of Theorem II in \cite{cerveau-berthier} (see pages 405,406 and 412-414).
Nevertheless, because of the extension problem to $X_f$, let us give details about the construction of the function $h$.

We take any point $ p \in X _f ^*= f^{-1}(0)\setminus \sing (X_f)$ and a small disk $\Sigma$, centered at $p$, transverse to $X_f$ (and therefore to $\omega$).


Because $X_f$ is irreducible, $X_f\setminus \sing (f)$ is connected and we conclude that, for $\Sigma$ small enough, the union
$\bigcup\limits_ {z \in \Sigma \setminus \{p\}} f^{-1}(f(z))$ is a neighborhood of the origin minus $X_f$, and \linebreak $\bigcup\limits_{z \in \Sigma} f^{-1}(f(z))$ is a neighborhood of the origin. Now we start by defining $h$ in $\Sigma$ as $h(p)=0$ and $h(z) = f(z), \forall z \in \Sigma$. We then extend $h$ to each fiber $f^{-1}(f(z)), z \ne p$ by integration, i.e.,
\[
h(w) = h(z) + \int \limits_z ^w \omega\big|_{f^{-1}(f(z))} = f(z) + \int\limits_z ^w\omega\big|_{f^{-1}(f(z))}, \, \forall w \in f^{-1}(f(z)).
\]

This line integral is well-defined due to the condition $\oint_\gamma \omega=0$ in (i) in the statement. Thus we have defined $h$ in $U\setminus X_f$ for some neighborhood $U$ of $X_f$ in $\mathbb C^n, 0$.

Notice that by definition we have $h$ holomorphic in the fibers $f^{-1}(z), z \ne 0$ and in the transverse disc $\Sigma$. So it is not difficult to conclude that (by a theorem of Hartogs) the function $h$ is holomorphic in $U\setminus X_f$. Now we observe that by construction $dh$ and $\omega$ coincide along the fibers $f^{-1}(z) , z \ne 0$. Therefore we can write $\omega = dh + a df$ for some holomorphic function $a \colon U \setminus X_f \to \mathbb C$ (notice that $df$ is non-singular in $U\setminus X_f$). Since $\omega \wedge df = dh \wedge df$ and $X_f$ is $\omega$-invariant, we conclude that
$(dh\wedge df) (z) \to 0$ as $ z \to X_f$. In other words $h(z) \to 0= h(p)$ as $z \to X_f$.
In particular $h\colon U \setminus X_f \to \mathbb C$ is bounded and by Riemann extension theorem $h$ admits an unique holomorphic extension to $X_f$. This extension satisfies $h(X_f)=\{0\}$.

Once we have $\omega = adf + dh$ with $\omega$ and $f, h$ holomorphic in $U$, the same holds for $a$ because $adf = \omega - dh$.

\end{proof}
This completes the proof of Proposition~\ref{Proposition:integration}.
\end{proof}

\begin{proof}[Proof of Lemma~\ref{Lemma:1}]
Lemma~\ref{Lemma:1} is now a straightforward consequence of Proposition~\ref{Proposition:integration} and L\^e-Saito Theorem (Theorem~\ref{Theorem:LeSaito}).
\end{proof}
\bibliographystyle{amsalpha}

\vglue.2in

\leftline{Dominique CERVEAU} \leftline{Universit\'e de Rennes I -
IRMAR} \leftline{Campus de Beaulieu - 35042 - Rennes Cedex}
\leftline{FRANCE}
\leftline{dominique.cerveau@univ-rennes1.fr}

\vglue.2in

\leftline{Bruno SC\'ARDUA} \leftline{Instituto de Matem\'atica}
\leftline{Universidade Federal do Rio de Janeiro} \leftline{Caixa
Postal 68530} \leftline{CEP. 21945-970 Rio de Janeiro - RJ}
\leftline{BRASIL}
\leftline{scardua@im.ufrj.br}

\end{document}